\documentclass[a4paper,10pt]{article}
\usepackage[T1]{fontenc}
\usepackage[latin1]{inputenc}
\usepackage{palatino, url, multicol}
\usepackage[sc]{mathpazo}
\usepackage{amsthm}
\usepackage{graphicx}
\usepackage{epsfig} 
\usepackage{amsmath}
\usepackage{latexsym}
\usepackage{amssymb}
\usepackage{amscd}
\usepackage{ifthen}
\usepackage{cite}
\newtheorem{theorem}{Theorem}[section]
\newtheorem{corollary}{Corollary}
\newtheorem{proposition}{Proposition}[section]

\theoremstyle{definition}

\theoremstyle{remark}
\newtheorem{remark}{Remark}
\title{Dynamics of the Tippe Top -- properties of numerical solutions versus the dynamical equations}

\author{
  Stefan Rauch-Wojciechowski, \hspace{8pt} Nils Rutstam \\
  \texttt{strau@mai.liu.se},\hspace{15pt}\texttt{nils.rutstam@liu.se}\\
  Matematiska institutionen, \\
  Link{\"o}pings universitet,\\
  SE--581 83 Link{\"o}ping, Sweden.\\
}

\begin{document}

\maketitle
{\abstract We study the relationship between numerical solutions for inverting Tippe Top and the structure of the dynamical equations. The numerical solutions confirm oscillatory behaviour of the inclination angle $\theta(t)$ for the symmetry axis of the Tippe Top. They also reveal further fine features of the dynamics of inverting solutions defining the time of inversion. These features are partially understood on the basis of the underlying dynamical equations.}
\bigskip

\noindent{\it Key words: tippe top; rigid body; nonholonomic mechanics; numerical solutions}

\section{Introduction}

\begin{figure}[ht]
\begin{center}
\includegraphics[scale=0.90]{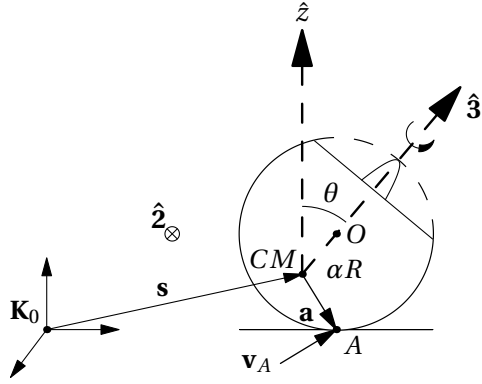}
\caption{Diagram of the TT. Note that $\mathbf{a}=R\alpha\mathbf{\hat{3}}-R\hat{z}$.\label{TT_diagram}}
\end{center}
\end{figure}

A Tippe Top (TT) is modeled by an axially symmetric sphere of mass $m$ and radius $R$ with center of mass ($CM$) shifted w.r.t. the geometrical center $O$ by $R\alpha$, $0<\alpha<1$ (see Fig. \ref{TT_diagram}). In a toy TT it is achieved by cutting off a slice of a sphere and by substituting it with a small peg that is used for spinning the TT. The sphere is rolling and gliding on a flat surface and subjected to the gravitational force $-mg\hat{z}$.

When spun slowly on the spherical part the TT spins wobbly for some time and comes to a standstill due to loss of energy caused by spinning friction. However when the initial spin is sufficiently fast the TT displays a counterintuitive behaviour of flipping upside down to spin on the peg with $CM$ above the geometrical center $O$. It continues spinning for some time until it falls down due to frictional loss of energy.
This flipping behaviour of TT we call {\it inversion}.

\bigskip

\noindent To describe motion of the TT we choose (as in \cite{Nisse2,Nisse3}) a fixed inertial reference frame $(\widehat{X},\widehat{Y},\widehat{Z})$ with $\widehat{X}$ and $\widehat{Y}$ parallel to the supporting plane and with vertical $\widehat{Z}$. We place the origin of this system in the supporting plane. Let $(\hat{x},\hat{y},\hat{z})$ be a frame defined through rotation around $\widehat{Z}$ by an angle $\varphi$, where $\varphi$ is the angle between the plane spanned by $\widehat{X}$ and $\widehat{Z}$ and the plane spanned by the points $CM$, $O$ and $A$.

The third reference frame $(\mathbf{\hat{1}},\mathbf{\hat{2}},\mathbf{\hat{3}})$, with origin at $CM$, is defined by rotating $(\hat{x},\hat{y},\hat{z})$ by an angle $\theta$ around $\hat{y}$. Thus $\mathbf{\hat{3}}$ is parallel to the symmetry axis, and $\theta$ is the angle between $\hat{z}$ and $\mathbf{\hat{3}}$. This frame is not fully fixed in the body. The axis $\mathbf{\hat{2}}$ points behind the plane of the picture of Fig.~\ref{TT_diagram}.   

We let $\mathbf{s}$ denote the position of $CM$ w.r.t. the origin of the frame $(\widehat{X},\widehat{Y},\widehat{Z})$ and the vector from $CM$ to $A$ is $\mathbf{a}=R(\alpha\mathbf{\hat{3}}-\hat{z})$.  
The orientation of the body w.r.t. the inertial reference frame $(\hat{X},\hat{Y},\hat{Z})$ is described by the Euler angles $(\theta,\varphi,\psi)$, where $\psi$ is the rotation angle of the sphere about the symmetry axis. With this notation, the angular velocity of the TT is $\boldsymbol{\omega}=-\dot{\varphi}\sin\theta\mathbf{\hat{1}}+\dot{\theta}\mathbf{\hat{2}}+(\dot{\psi}+\dot{\varphi}\cos\theta)\mathbf{\hat{3}}$, and we denote $\omega_3 :=\dot{\psi}+\dot{\varphi}\cos\theta$. 

The principal moments of inertia along the axes $(\mathbf{\hat{1}},\mathbf{\hat{2}},\mathbf{\hat{3}})$ are denoted by $I_1=I_2$ and $I_3$, so the inertia tensor $\mathbb{I}$ will have components $(I_1,I_1,I_3)$ with respect to the $(\mathbf{\hat{1}},\mathbf{\hat{2}},\mathbf{\hat{3}})$-frame. The axes $\mathbf{\hat{1}}$ and $\mathbf{\hat{2}}$ will be principal axes due to the axial symmetry of TT.

Motion of TT is described by the standard Newton equations for a rolling and gliding rigid body. They have the vector form 
\begin{equation}
\label{TTequ} m\mathbf{\ddot{s}}=\mathbf{F}-mg\hat{z},\quad \mathbf{\dot{L}}=\mathbf{a}\times\mathbf{F},\quad\mathbf{\dot{\hat{3}}}=\boldsymbol{\omega}\times\mathbf{\hat{3}}=\frac{1}{I_1}\left(\mathbf{L}\times\mathbf{\hat{3}}\right),
\end{equation} 
where $\mathbf{F}=\mathbf{F}_{R}+\mathbf{F}_{f}=g_n\hat{z}-\mu g_n\mathbf{v}_A$ is the external force acting on the TT at the point of support $A$. In this model $\mathbf{F}=\mathbf{F}_{R}+\mathbf{F}_{f}$ consists of a vertical normal force $g_n\geq 0$ and a viscous-type friction force $-\mu g_{n}\mathbf{v}_A$, acting against the gliding velocity $\mathbf{v}_A$, with $\mu\geq 0$ the friction coefficient.

Equations \eqref{TTequ} admit Jellett's integral of motion $\lambda=-\mathbf{L\cdot a}$, $\dot{\lambda}=-\mathbf{\dot{L}}\cdot\mathbf{a}-\mathbf{L}\cdot\mathbf{\dot{a}}=-(\mathbf{a\times F})\cdot\mathbf{a}-\mathbf{L}\cdot(\frac{R\alpha}{I_1}(\mathbf{L}\times\mathbf{\hat{3}}))=0$. The total energy $E=\frac{1}{2}m\mathbf{\dot{s}}^2+\frac{1}{2}\boldsymbol{\omega}\cdot\mathbf{L}+mg\mathbf{s}\cdot\hat{z}$ is monotonically decreasing since $\dot{E}=\mathbf{v}_{A}\cdot\mathbf{F}=-\mu g_n|\mathbf{v}_A|^2$. The rolling and gliding solutions satisfy the one-sided constraint $(\mathbf{s+a})\cdot\hat{z}=0$ and its derivative satisfy $\mathbf{v}_{A}\cdot\hat{z}=0$. When Eqs. \eqref{TTequ} are expressed in terms of the Euler angles we get \cite{Nisse2,Nisse4}:
\begin{align}
\label{ddth}\ddot{\theta}=&\frac{\sin\theta}{I_1}\left(I_1\dot{\varphi}^2\cos\theta-I_3\omega_3\dot{\varphi}-R\alpha g_n\right)+\frac{R\mu g_n\nu_x}{I_1}(1-\alpha\cos\theta),&\\
\label{ddph}\ddot{\varphi}=&\frac{I_3\dot{\theta}\omega_3-2I_1\dot{\theta}\dot{\varphi}\cos\theta-\mu g_n\nu_y R(\alpha-\cos\theta)}{I_1\sin\theta},&\\
\label{dot_om}\dot{\omega}_3=&-\frac{\mu g_n\nu_y R\sin\theta}{I_3},&\\
\label{nu_x}\dot{\nu}_x=&\frac{R\sin\theta}{I_1}\left(\dot{\varphi}\omega_3\left(I_3(1-\alpha\cos\theta)-I_1\right)+g_nR\alpha(1-\alpha\cos\theta)-I_1\alpha(\dot{\theta}^2+\dot{\varphi}^2\sin^2\theta)\right)\nonumber\\
&\qquad-\frac{\mu g_n\nu_x}{mI_1}\left( I_1+mR^2(1-\alpha\cos\theta)^2\right)+\dot{\varphi}\nu_y,\\
\label{nu_y}\dot{\nu}_y=&-\frac{\mu g_n\nu_y}{mI_1 I_3}\left(I_1I_3+mR^2 I_3(\alpha-\cos\theta)^2+mR^2I_1\sin^2\theta\right)\nonumber\\ 
&+\frac{\omega_3\dot{\theta} R}{I_1}\left(I_3(\alpha-\cos\theta)+I_1\cos\theta\right)-\dot{\varphi}\nu_x.
\end{align}
These equations are a complicated nonlinear dynamical system for 6 unknowns.  The value of the normal force $g_n$ in the dynamical equations is determined by the second derivative $\frac{d^2}{dt^2}(\mathbf{s+a})\cdot\hat{z}=0$ of the contact constraint:
\begin{equation}
\label{equ_for_g_n}g_n=\frac{mgI_1+mR\alpha(\cos\theta(I_1\dot{\varphi}^2\sin^2\theta+I_1\dot{\theta}^2)-I_3\dot{\varphi}\omega_3\sin^2\theta)}{I_1+mR^2\alpha^2\sin^2\theta-mR^2\alpha\sin\theta(1-\alpha\cos\theta)\mu\nu_x}.
\end{equation}
Equations \eqref{ddth}--\eqref{nu_y} admit the same Jellett integral that expressed through Euler angles reads $\lambda=RI_1\dot{\varphi}\sin^2\theta-I_3\omega_3(\alpha-\cos\theta)$.

\bigskip

\noindent An asymptotic analysis of TT inversion \cite{Eben,Mars,RSG} provides sufficient conditions for physical parameters that the TT have to satisfy and for initial conditions so that TT is inverting.
\begin{theorem}\label{inversion_cond}
For a Tippe Top with parameters satisfying $1-\alpha<\gamma=\frac{I_1}{I_3}<1+\alpha$, an inverted spinning solution is the only Lyapunov stable state in the asymptotic set of nongliding solutions $M=\{(\mathbf{L},\mathbf{\hat{3}},\mathbf{v}_A):\dot{E}=-\mu g_n|\mathbf{v}_A|^2=0\}$ provided that $\lambda>\max\left\{\lambda_{\text{thres}}=\frac{\sqrt{mgR^3\alpha I_3}(1+\alpha)^2}{\sqrt{1+\alpha-\gamma}},\lambda_{\text{up}}=\frac{\sqrt{mgR^3\alpha I_3}(1-\alpha)^2}{\sqrt{\alpha+\gamma-1}}\right\}$.
\end{theorem}

\begin{remark}
When $1-\alpha<1-\alpha^2<\gamma<1+\alpha$, then $\lambda_{\text{thres}}>\lambda_{\text{up}}$.
\end{remark}
\begin{remark}
Inversion of TT for $\lambda>\max\{\lambda_{\text{thres}},\lambda_{\text{up}}\}$ is a consequence of the LaSalle theorem applied to each solution having positive value $g_{n}(t)>0$ of the normal force \cite{RSG}. A direct application of the LaSalle theorem would require specification of an initial compact invariant set, which is difficult to define due to the one-sided constraint $(\mathbf{s+a})\cdot\hat{z}=0$ that does not exclude existence of solutions having negative $g_n(t)<0$.
\end{remark}

\begin{corollary}
When $\lambda>\max\{\lambda_{\text{thres}},\lambda_{\text{up}}\}$ then the asymptotic set $M$ contains only one solution and every solution with $g_n(t)>0$ is asymptotically approaching the inverted spinning solution $(\mathbf{L},\mathbf{\hat{3}},\mathbf{v}_A)=\left(\frac{\lambda}{R(1+\alpha)}\hat{z},\hat{z},0\right)$.
\end{corollary}

Application of the LaSalle theorem to TT inversion provides only an existential result. It states when a TT inverts but it says nothing about dynamical behaviour of solutions during the inversion.

There have been attempts to study dynamics of inversion by applying a gyroscopic balance condition \cite{Moff,Ued} where the quantity $\xi=I_3\omega_3-I_1\dot{\varphi}\cos\theta$ is assumed to be close to zero. As has been pointed out in these articles, this condition is useful for explaining rising of a rotating egg. For TT it is approximately satisfied in a certain neighborhood of the angle $\theta=\frac{\pi}{2}$ when the equator of the TT is in contact with the supporting plane. The use of TT equations simplified through the condition $\xi\approx0$ leads to an oversimplified equation of the form $\dot{\theta}=\frac{\mu g_n R^2(1-\alpha\cos\theta)^2}{\lambda}\nu_y$  providing a monotonously increasing $\theta(t)$ as a solution \cite{Moff}. It is intuitively comprehensible that such a solution may reflect correctly some type of averaged behaviour of the inclination angle $\theta(t)$ during inversion. The problem is however that no suitable definition of an averaged angle $\theta(t)$ is available. As is well known from numerical simulations $\theta(t)$ oscillates (meaning that $\dot{\theta}(t)$ changes sign many times) about a logistic type curve during inversion. 

\bigskip

\noindent A rigorous approach showing oscillatory behaviour of TT inverting solutions has been proposed in \cite{Rau,Nisse2,Nisse4}. It is based on an integrated form of TT equations that are equivalent to the Euler angle equations and to the dynamical Eqs. \eqref{TTequ}. We get the integrated form of TT equations by considering functions which are integrals of motion for the rolling axisymmetric sphere, the modified energy (with $\mathbf{\dot{s}}=-\boldsymbol{\omega}\times\mathbf{a}$): 
\begin{align}
\label{Etilde_def}\tilde{E}=&\frac{1}{2}m\mathbf{\dot{s}}^2+\frac{1}{2}\boldsymbol{\omega}\cdot\mathbf{L}+mg\mathbf{s}\cdot\hat{z}=\frac{1}{2}\left(I_1\dot{\varphi}^2\sin^2\theta+I_1\dot{\theta}^2+I_3\omega_{3}^2\right)+mgR(1-\alpha\cos\theta)\nonumber\\
&+\frac{1}{2}mR^2\bigg[(\alpha-\cos\theta)^2(\dot{\theta}^2+\dot{\varphi}^2\sin^2\theta)+\sin^2\theta(\dot{\theta}^2+\omega_{3}^2+2\omega_3\dot{\varphi}(\alpha-\cos\theta))\bigg].
\end{align} 
and the Routh function:
\begin{equation}
\label{Routh_def}D(\theta,\omega_3)=\omega_3\sqrt{I_1I_3+mR^2I_3(\alpha-\cos\theta)^2+mR^2I_1\sin^2\theta}=I_3\omega_3\sqrt{d(\cos\theta)},
\end{equation} 
where $d(\cos\theta)=\gamma+\sigma(\alpha-\cos\theta)^2+\sigma\gamma(1-\cos^2\theta)$, $\sigma=\frac{mR^2}{I_3}$.
For the rolling and gliding TT, equations of motion \eqref{ddth}--\eqref{nu_y} are equivalent to the following equations:
\begin{align}
&\frac{d}{dt}\lambda(\theta,\dot{\theta},\dot{\varphi},\omega_3)=0,\\
&\frac{d}{dt}D(\theta,\omega_3)=\frac{\gamma m}{\alpha\sqrt{d(\hat{z}\cdot\mathbf{\hat{3}})}}(\hat{z}\times\mathbf{a})\cdot\mathbf{\dot{v}}_A,\\
&\frac{d}{dt}\tilde{E}(\theta,\dot{\theta},\dot{\varphi},\omega_3)=m(\boldsymbol{\omega}\times\mathbf{a})\cdot\mathbf{\dot{v}}_A\\
&\frac{d}{dt}m\mathbf{\dot{r}}=-\mu g_n\mathbf{v}_A,
\end{align}
where $\mathbf{r}=\mathbf{s}-s_{\hat{z}}\hat{z}$. From these equations it follows that for any given solution $(\theta(t),\dot{\theta}(t),\dot{\varphi}(t),\omega_3(t),\nu_x(t),\nu_y(t))$ the functions $D(t)$, $\tilde{E}(t)$ depend on time as $D(t)=D(\theta(t),\omega_3(t))$, $\tilde{E}(t)=\tilde{E}(\theta(t),\dot{\theta}(t),\dot{\varphi}(t),\omega_3)$ and then {\it this} solution satisfies the three equations
\begin{align}
\label{lambda_equ}\lambda &=RI_1\dot{\varphi}\sin^2\theta-RI_3\omega_3(\alpha-\cos\theta),\\
\label{Routh_equ}D(t)&=I_3\omega_3\sqrt{d(\cos\theta)},\\
\label{E_equ}\tilde{E}(t)&=\frac{1}{2}\left(I_1\dot{\varphi}^2\sin^2\theta+I_1\dot{\theta}^2+I_3\omega_{3}^2\right)+mgR(1-\alpha\cos\theta)\nonumber\\
&+\frac{1}{2}mR^2\bigg[(\alpha-\cos\theta)^2(\dot{\theta}^2+\dot{\varphi}^2\sin^2\theta)+\sin^2\theta(\dot{\theta}^2+\omega_{3}^2+2\omega_3\dot{\varphi}(\alpha-\cos\theta))\bigg],
\end{align}
and the equation $m\mathbf{\ddot{r}}=-\mu g_n\mathbf{v}_A$.
The three first equations are decoupled from the last one and on their own they recall the problem of a purely rolling axisymmetric sphere solved by Chaplygin \cite{Chap}. By eliminating $\dot{\varphi}(t)$, $\omega_3(t)$ from \eqref{lambda_equ} and \eqref{Routh_equ} and substituting into \eqref{E_equ} one obtains a single first order differential equation for $\theta(t)$:
\begin{equation}
\label{METT}\tilde{E}(t)=g(\cos\theta)\dot{\theta}^2+V(\cos\theta,D(t),\lambda),
\end{equation}
where $g(\cos\theta)=\frac{1}{2}I_3\left(\sigma((\alpha-\cos\theta)^2+1-\cos^2\theta)+\gamma\right)$ and
\begin{equation}
V(z,D(t),\lambda)=mgR(1-\alpha z)+\frac{(\lambda\sqrt{d(z)}+RD(t)(\alpha-z))^2}{2I_3R^2\gamma^2(1-z^2)}+\frac{(R^2D(t)^2-\sigma\lambda^2)}{2R^2I_1}.
\end{equation}
The equation \eqref{METT} has the same algebraic form as the separation equation for a rolling sphere but here $D(t)$, $\tilde{E}(t)$ are time-dependent, so the equation is not separable or explicitly solvable. It is called the {\it Main Equation for the Tippe Top}. 

Each solution of the TT equations satisfies its own main equation with suitable functions $D(t)$, $\tilde{E}(t)$ which are {\it a priori unknown}. But all these equations have the same algebraic form and they can be understood as describing a particle with variable mass $2g(\cos\theta)$ moving in a potential well $V(\cos\theta,D(t),\lambda)$ that is deforming in time. In a generic situation the particle reflects many times between the walls of $V$ and, due to energy dissipation, its position $\theta(t)$ goes toward the minimum $\theta_{\min}(t)$ of the potential.

By the LaSalle theorem every Tippe Top satisfying theorem \ref{inversion_cond} {\it has} to invert. For inverting solutions $\lim_{t\to-\infty}\theta(t)=0$, $\lim_{t\to\infty}\theta(t)=\pi$, $\lim_{t\to-\infty}\mathbf{L}(t)=L_0\hat{z}$ and $\lim_{t\to\infty}\mathbf{L}(t)=L_1\hat{z}$ and all pairs $D(t)$, $\tilde{E}(t)$ satisfy the same asymptotic conditions that can be found from the Jellett integral: $\lambda=R(1-\alpha)L_0=R(1+\alpha)L_1$. So $\lim_{t\to-\infty}D(t)=D_0=\frac{\lambda\sqrt{d(1)}}{R(1-\alpha)}$, $\lim_{t\to\infty}D(t)=D_1=-\frac{\lambda\sqrt{d(-1)}}{R(1+\alpha)}$ and $\lim_{t\to-\infty}\tilde{E}(t)=\tilde{E}_0=\frac{\lambda^2}{2R^2I_3(1-\alpha)^2}+mgR(1-\alpha)$, $\lim_{t\to\infty}\tilde{E}(t)=\tilde{E}_1=\frac{\lambda^2}{2R^2I_3(1+\alpha)^2}+mgR(1+\alpha)$.

This means that in the plane $(D,\tilde{E})$ each inverting solution draws a curve $(D(t),\tilde{E}(t))$ that is moving in finite time from a small neighborhood of $(D_0,\tilde{E}_0)$ to a small neighborhood of $(D_1,\tilde{E}_1)$, as the effective potential $V(\cos\theta,D(t),\lambda)$ is deformed by $D(t)$.

\bigskip
 
\noindent In \cite{Nisse4} the deformation is considered in the special case when the effective potential $V(z,D(t),\lambda)$ is a rational function of $\cos\theta$. This occurs when the physical parameters satisfy $1-\alpha^2<\gamma<1$, $\sigma=\frac{1-\alpha}{\gamma+\alpha^2-1}$, conditions that may be realized by parameters of a real toy TT. The analysis of the potential is then simplified, and it has been proved that
\begin{proposition}\label{nbhds_of_poles}
Assume that $\lambda>\lambda_{\text{thres}}=\frac{\sqrt{mgR^3I_3\alpha}(1+\alpha)^2}{\sqrt{1+\alpha-\gamma}}$.
\begin{itemize}
\item[i)] For any (small) $\epsilon>0$ there is a $\delta_{-}(\epsilon,\lambda)>0$ such that for every positive $\delta<\delta_{-}(\epsilon,\lambda)$ the potential $V(z,D,\lambda)$ has a minimum $z_{\min}$ in the interval $[-1,-1+\epsilon]$ for $D=D_1+\frac{\delta}{R(1+\alpha)\sqrt{\gamma+\alpha^2-1}}$  . 

\item[ii)] For any (small) $\epsilon>0$ there is a $\delta_{+}(\epsilon,\lambda)>0$ such that for every positive $\delta<\delta_{+}(\epsilon,\lambda)$ the potential $V(z,D,\lambda)$ has a minimum $z_{\min}$ in the interval $[1-\epsilon,1]$ for $D=D_0-\frac{\delta}{R(1-\alpha)\sqrt{\gamma+\alpha^2-1}}$.
\end{itemize}
\end{proposition} 
This proposition states that if $\lambda>\lambda_{\text{thres}}$, so that TT is inverting, then the minimum $z_{\min}$ of $V(z,D(t),\lambda)$ is moving from an $\epsilon$-neighborhood of $z=1$ to an $\epsilon$-neighborhood of $z=-1$. In \cite{Nisse4} it has also been shown that the time of passage for $\theta(t)$ between two turning points has a uniform bound $T_{\text{upp}}=21.95\left(\frac{RI_3\gamma(\alpha+1-\gamma)}{\alpha\lambda+\alpha RD\sqrt{\gamma+\alpha^2-1}}\right)$. If the time of inversion is an order of magnitude larger $T_{\text{inv}}>10 T_{\text{upp}}$, then the angle $\theta(t)$ has to perform several oscillations during the inversion. 
The motion of the symmetry axis $\mathbf{\hat{3}}(t)$ can be seen as nutational motion within a belt that is moving from the neighborhood of the north pole of the unit sphere $S^2$ to the neighborhood of the south pole. As an inverting solution is approaching the inverted asymptotic state $(\mathbf{L},\mathbf{\hat{3}},\mathbf{v}_A)=\left(\frac{\lambda}{R(1+\alpha)}\hat{z},\hat{z},0\right)$ the velocity $\dot{\theta}(t)\to 0$ and the total energy $E(t)\to E_1=\lim_{t\to\infty}V(z_{\min}(t),D(t),\lambda)$.

This qualitative analysis of inverting solutions shows oscillatory character of inverting solutions and no more. The dynamical equations make it difficult to derive rigorous statements about further properties of solutions.

\bigskip

\noindent As we shall show below these properties of inverting solutions are well confirmed by numerical simulations made with the use of the Python 2.7 open source library SciPy \cite{Scipy} for realistic values of parameters corresponding to a typical toy TT available commercially. These simulations also show that generic dynamics of inverting TT has other distinctive features that does not follow from the theoretical analysis presented above. They have yet to be found from the underlying dynamical equations and such derivation is a daunting challenge.

In order to understand better these fine features of inverting solutions we adopt a converse approach of analysing how the results of simulations {\it agree} with the dynamical equations.

It should be stated clearly that such analysis does not provide proofs of existence of these new features but it improves understanding of relationships that hold between dynamical variables of an inverting TT. This enables formulation of some specific hypotheses that may have a chance of being rigorously proved from the dynamical equations.

\section{Simulations of a toy TT}

\subsection{Behaviour of a reference inverting solution}

We have taken realistic values of the physical parameters for a typical toy TT; $m=0.02$ kg, $R=0.02$ m, $\alpha=0.3$, $I_3=\frac{2}{5}mR^2$, $I_1=\frac{131}{350}mR^2$ and $g=9.82\text{ m/s}^2$. They satisfy the rationality condition $\frac{mR^2}{I_3}=\frac{1-\alpha}{\gamma+\alpha^2-1}$. We also choose $\mu=0.3$. These values are close to the parameter values used in simulations by other authors \cite{Coh,Ued} . 
We study initial conditions (IC) with $\lambda\approx 2\lambda_{\text{thres}}$ and with small $\theta(0)=0.1$ rad that give inversion of the TT. As reference initial conditions we take (similarly as in \cite{Coh,Ued})
\begin{equation}
\label{IC_ref}IC_{\text{ref}}(0)=\{\theta(0)=0.1 \text{ rad},\;\dot{\theta}(0)=0,\;\dot{\varphi}(0)=0,\;\omega_3(0)=155 \text{ rad}/s,\;\nu_x(0)=0,\;\nu_y(0)=0\}
\end{equation}
The use of IC with $\dot{\varphi}(0)\approx0$ is justified by a typical demonstration of a toy TT that is initially started with strong spin (here $\dot{\psi}(0)=155$ rad/s) about the axis $\mathbf{\hat{3}}$ before it hits the supporting plane.

The consistency of calculated values of the dynamical variables has been checked by finding that the value of Jellett's integral $\lambda=RI_1\dot{\varphi}\sin^2\theta-RI_3\omega_3(\alpha-\cos\theta)$ is constant with precision of order $10^{-6}\lambda$. The condition $\lambda=\text{constant}$ essentially connects $\dot{\varphi}$ and $\omega_3=\dot{\psi}+\dot{\varphi}\cos\theta$, and initially $\omega_3(0)\approx\dot{\varphi}(0)+\dot{\psi}(0)$ since $\theta(0)=0.1$ and $\cos\theta\approx 1$.

\begin{figure}[tbp]
\begin{center}
\includegraphics[width=\textwidth,height=0.6\textheight]{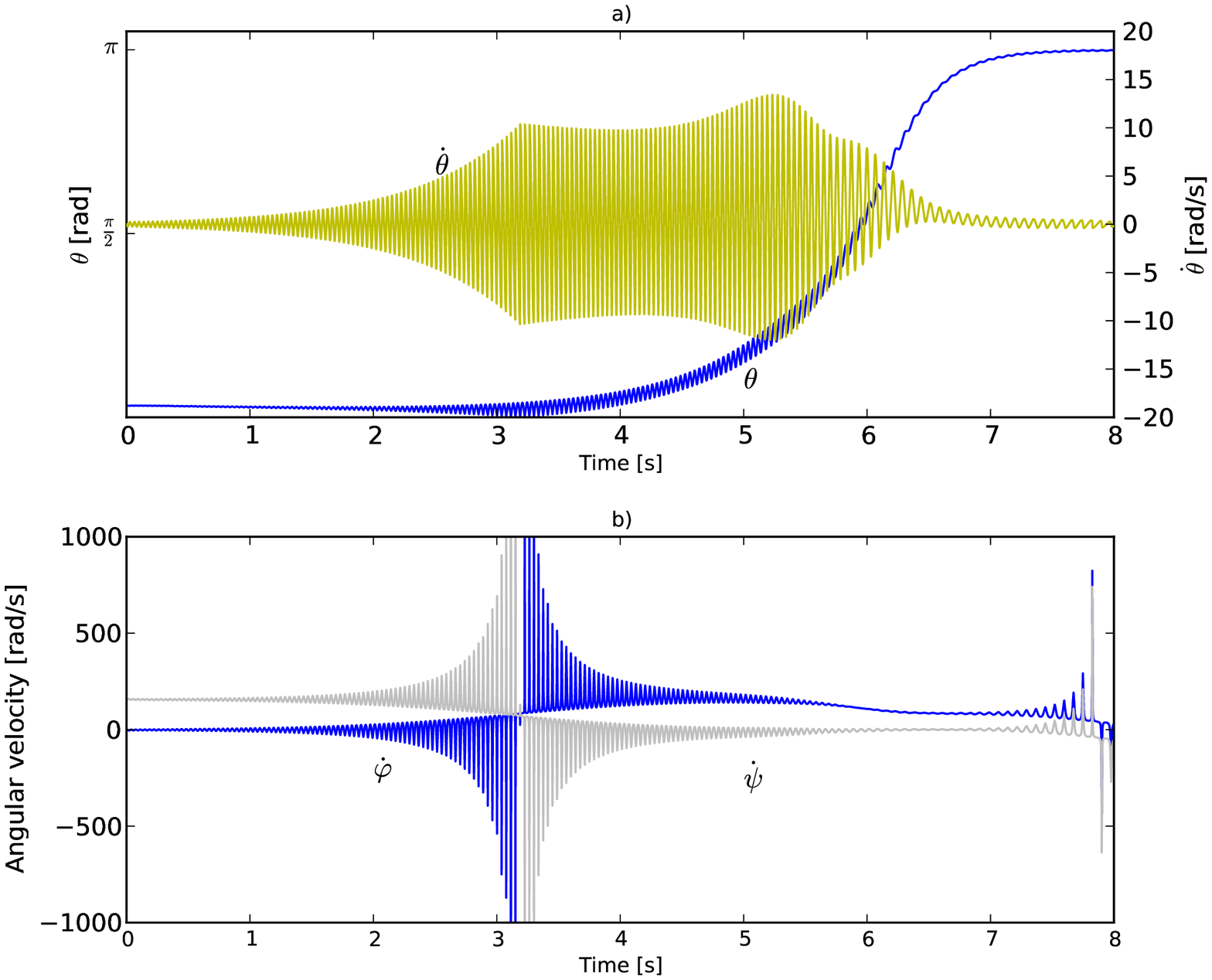}
\caption{Plots obtained by numerically integrating Eqs. \eqref{ddth}-\eqref{nu_y} and Eq. \eqref{equ_for_g_n} using the Python 2.7 open source library SciPy \cite{Scipy}. Values of physical parameters are taken as in section 2.1. Initial values are $\theta(0)=0.1$ rad, $\omega_3(0)=155.0$ rad/s, $\dot{\varphi}(0)=\dot{\theta}(0)=0$ and $\nu_x(0)=\nu_y(0)=0$. Plot $a$ shows the evolution of the inclination angle $\theta(t)$ (blue) and the angular velocity $\dot{\theta}(t)$ (yellow) and plot $b$ show the evolution of the angular velocities $\dot{\varphi}(t)$ (blue) and $\dot{\psi}(t)$ (grey).\label{total_plots1}}
\end{center}
\end{figure}

Quite remarkably inverting solutions with $\lambda\approx 2\lambda$ and small $\dot{\varphi}(0)\approx0$ display the same universal behaviour (see Fig. \ref{total_plots1}) consisting of three main phases separated by two moments $t_{\theta\approx0}$, $t_{\theta\approx\pi}$ of singular behaviour of $\dot{\varphi}(t)$ when this angular velocity changes sign from large negative values to large positive values. At the same moments the angle $\theta(t)$ takes correspondingly locally minimal value close to 0 at $t_{\theta\approx0}$ and then locally maximal value close to $\pi$ at $t_{\theta\approx\pi}$. 

\bigskip

\noindent The three phases are thus: 
\begin{itemize}
\item[1)] Initial sychronisation of dynamical variables (that ends with singular values for angular velocities $\dot{\varphi}$, $\dot{\psi}$).
\item[2)] Climbing phase, when $\theta(t)$ climbs in an oscillatory way from $\theta\approx 0$ towards values of $\theta$ close to $\theta\approx \pi$. It ends with singular values of angular velocities $\dot{\varphi}$, $\dot{\psi}$.
\item[3)] Asymptotic stabilisation phase when $\theta(t)$ slowly approaches $\theta\approx\pi$.
\end{itemize}

The existence of singular behaviour is not affected by small perturbations of IC. Thus $t_{\theta\approx0}$ may be called an {\it initiation time} and $t_{\theta\approx\pi}$ an {\it ending time} for the climbing of $\theta(t)$. A {\it climbing time} is the quantity $T_{\text{inv}}=t_{\theta\approx \pi}-t_{\theta\approx0}$ from the first definite change of sign of $\dot{\varphi}(t)$ at $\theta(t)\approx0$ to the second change of sign of $\dot{\varphi}(t)$ at $\theta(t)\approx\pi$.

The initial synchronisation phase is the time required to accelerate the angular velocity $\dot{\varphi}$ to about $80$ rad/s needed for the climbing phase to start.
When values of dynamical variables attained close to the first time of singularity are taken as initial conditions
\begin{equation}
IC(0)=IC(t_{\theta\approx0})=\{\theta(t_{\theta\approx 0}),\;\dot{\theta}(t_{\theta\approx0}),\;\dot{\varphi}(t_{\theta\approx0}),\;\omega_3(t_{\theta\approx0}),\;\nu_x(t_{\theta\approx0}),\;\nu_y(t_{\theta\approx0})\},
\end{equation}
the solution $\theta(t)$ climbs immediately.

The graph for the inclination angle $\theta(t)$ in Fig. \ref{total_plots1}a displays the main features of an inverting TT. It has a general form of a logistic type curve superposed with small amplitude oscillations. These oscillations are visible in the graph of $\dot{\theta}(t)$ that changes sign frequently. Zeros of $\dot{\theta}(t)$ correspond to subsequent local maximum and minimum points of the $\theta(t)$-curve (see Fig. \ref{total_plots1}a). These properties confirm the picture implied by the Main Equation for the TT saying that the symmetry axis $\mathbf{\hat{3}}(t)$ performs nutational motion within a narrow band that, during the inversion, is shifting on the unit sphere $S^2$ from a neighborhood of the north pole to a neighborhood of the south pole. Large values of $\dot{\theta}(t)$ during the climbing phase reflect widening of the nutational band when $\theta(t)$ is crossing $\frac{\pi}{2}$. The band is narrowing (and $\dot{\theta}(t)\to 0$) when $\mathbf{\hat{3}}(t)$ is approaching the south pole of $S^2$.   
 
In the graph \ref{total_plots1}b for $\dot{\varphi}(t)$ and $\dot{\psi}(t)$ we can discern an unexpected phenomenon (also visible in corresponding graphs published in \cite{Coh,Ued}) of sudden increase, by orders of magnitude, of oscillation amplitudes in vicinity of the initiation time $t_{\theta\approx0}\approx 3.2$ s and in vicinity of the ending time $t_{\theta\approx\pi}\approx 7.87$ s. Both times are clearly distinguished by a change of sign of high amplitude oscillations.

\begin{figure}[tbp]
\begin{center}
\includegraphics[width=\textwidth,height=0.6\textheight]{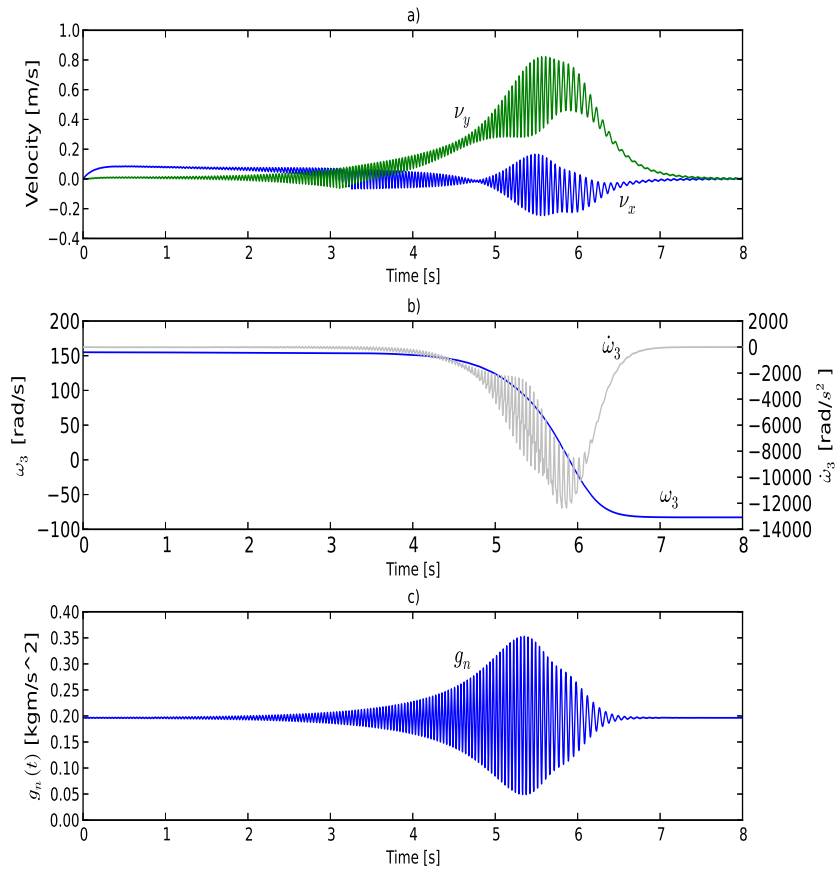}
\caption{Plots obtained as in figure \ref{total_plots1}. Plot $a$ shows the gliding velocities $\nu_x(t)$ (blue) and $\nu_y(t)$ (green), plot $b$ shows the evolution of the angular velocity $\omega_3(t)$ (blue) and the derivative $\dot{\omega}_3(t)$ (grey) and plot $c$ show the evolution of the value of the normal force $g_n\hat{z}$ which is positive at all times and oscillates about $mg=0.2$.\label{total_plots2}}
\end{center}
\end{figure}

The strong oscillations of $\dot{\varphi}(t)$, $\dot{\psi}(t)$ in vicinity of $t_{\theta\approx0}$ and of $t_{\theta\approx\pi}$, visible in Fig. \ref{total_plots1}b, compensate each other because the variable $\omega_3(t)=\dot{\psi}(t)+\dot{\varphi}(t)\cos\theta(t)$ (Fig. \ref{total_plots2}b) is (unexpectedly) an almost monotonously decreasing function from the initial value of $155$ rad/s to about $-85$ rad/s. 
The decreasing behaviour of $\omega_3$ reflects inversion of direction of axis $\mathbf{\hat{3}}$ (with respect to the direction of $\mathbf{L}$ and $\boldsymbol{\omega}$) and also reflects frictional loss of rotational energy. The curve $\omega_3(t)$ is in the beginning and at the end almost horizontal. During the climbing phase of TT between 3.2 and 7.86 seconds the symmetry axis $\mathbf{\hat{3}}$ turns upside down, so the projection $\boldsymbol{\omega}\cdot\mathbf{\hat{3}}$ of the predominantly vertical angular velocity $\boldsymbol{\omega}$ on $\mathbf{\hat{3}}$ changes sign. 

The function $\omega_3(t)$ is not, however, everywhere monotonously decreasing because its derivative $\dot{\omega}_3=-\frac{R\mu g_n}{I_3}\nu_y\sin\theta$ may take positive values when $\nu_y$ becomes negative in a neighborhood of $t_{\theta\approx 0}$. Indeed, the graph for $\dot{\omega}_3(t)$ in Fig. \ref{total_plots2}b acquires positive spikes in a neighborhood of $t_{\theta\approx 0}$.

The value $g_n(t)$ of the reaction force $g_n(t)\hat{z}$ stays positive all time during the inversion (Fig. \ref{total_plots2}c) but it oscillates strongly during the climbing phase. So these solutions fulfill the assumption needed for validity of the LaSalle type theorem \ref{inversion_cond}. Values of $g_n(t)$ are oscillating about $mg$ -- the value of the reaction force when a static TT is standing on its bottom.

\bigskip

\noindent Further information about behaviour of $\dot{\varphi}(t)$ can be deduced from the expression for Jellett's integral \eqref{lambda_equ} rewritten as
\begin{align}
\label{phi_dot_equ}\dot{\varphi}(t)\sin^2\theta(t)&=\frac{\lambda}{RI_3}+\omega_3(t)(\alpha-\cos\theta(t))=-\omega_3(0)(\alpha-\cos\theta(0))+\omega_3(t)(\alpha-\cos\theta(t))\nonumber\\
&=\omega_3(0)(\alpha-\cos\theta(0))\left(1-\frac{\omega_3(t)(\alpha-\cos\theta(t))}{\omega_3(0)(\alpha-\cos\theta(0))}\right),
\end{align}
since $\lambda=-RI_3\omega_3(0)(\alpha-\cos\theta(0))$. In the vicinity of $t_{\theta\approx 0}$ the r.h.s. of \eqref{phi_dot_equ} is a small quantity because in the quotient $\omega_3(t)$ ($\approx 154$ rad/s) is close to $\omega_3(0)=155$ rad/s and both angles $\theta(t)$, $\theta(0)$ are small.

The formula \eqref{phi_dot_equ} makes it easier to understand the interplay between $\theta(t)$ and $\dot{\varphi}(t)$ during the inversion. Since $\sin^2\theta\geq 0$ the sign of $\dot{\varphi}(t)$ is the same as the sign of the r.h.s. of \eqref{phi_dot_equ}. The part of the $\dot{\varphi}(t)$-graph left of $t_{\theta\approx 0}\approx 3.2$s shows that the quantity $\left(1-\frac{\omega_3(t)(\alpha-\cos\theta(t))}{\omega_3(0)(\alpha-\cos\theta(0))}\right)$ changes sign many times. As the r.h.s. of \eqref{phi_dot_equ} is small, because the quotient is close to 1, the amplitude of $\dot{\varphi}(t)$ can increase on approaching $t_{\theta\approx 0}$ only if $\sin^2\theta(t)$ becomes very small, and the axis $\mathbf{\hat{3}}$ almost hits the north pole of $S^2$. At time $t_{\theta\approx 0}$ the term $\left(1-\frac{\omega_3(t)(\alpha-\cos\theta(t))}{\omega_3(0)(\alpha-\cos\theta(0))}\right)$ acquires a positive sign and the amplitude of $\dot{\varphi}(t)$ decreases with rising angle $\theta(t)$. This is the beginning of inversion.

The moment of initiation of inversion of TT at $3.2$s noticed at graphs for $\dot{\varphi}(t)$ and $\dot{\psi}(t)$ is also visible in the graph \ref{total_plots2}a for the gliding velocities $\nu_x(t)$, $\nu_{y}(t)$. Both velocities are oscillatory but initially $\nu_{x}$ (blue) is positive and has larger value than $\nu_y$ (green) that oscillates close to $0$. 

Remarkably the graph of $\nu_y(t)$ crosses the graph of $\nu_x(t)$ at $t_{\theta\approx 0}\approx 3.2$s and from this moment the minimal of $\nu_y(t)$ are consistently higher than the minima of $\nu_x(t)$. During the climbing phase the oscillatory function $\nu_y(t)$ is positive and its values increase by about two orders of magnitude in comparison with values during the synchronatisation phase. Oscillations of $\nu_x(t)$ also increase and the mean value of $\nu_x$ becomes slightly negative. When the climbing is finished at $t_{\theta\approx\pi}\approx 7.86$s the graph of $\nu_y$ again crosses the graph of $\nu_x$ from above and minima of $\nu_y(t)$ become consistently lower than minima of $\nu_x(t)$.

\subsection{Testing initial conditions for inverting solutions}

The reference solution with IC \eqref{IC_ref} corresponds to a TT spun rapidly about the $\mathbf{\hat{3}}$ axis with $\omega_3(0)=\dot{\psi}(0)=155$ rad/s and launched upon a table with a small inclination angle $\theta(0)=0.1$ rad. We vary IC while keeping the value of Jellett's integral $\lambda\approx 2\lambda_{\text{thres}}$ above the threshold value, which is sufficient for inversion of TT.

As $\lambda=RI_1\dot{\varphi}\sin^2\theta-RI_3\omega_3(\alpha-\cos\theta)$ (with $\omega_3=\dot{\psi}+\dot{\varphi}\cos\theta$) does not depend on $\dot{\theta}$, $\nu_x$ or $\nu_y$ it is natural to consider first how, for fixed $\theta(0)=0.1$ rad, the initial distribution of the angular velocity between $\omega_3(0)$ and $\dot{\varphi}(0)$ affects the character of inverting solutions and then to study additional influence of nonvanishing IC for $\nu_x(0)$, $\nu_y(0)$ and $\dot{\theta}(0)$. Changes of the initial inclination angle $\theta(0)$ are also tested.

For testing influence of nonzero initial angular velocity $\dot{\varphi}(0)$ we consider $\dot{\varphi}(0)$ being maximally of the same order of magnitude as $\omega_3(0)$ and belonging to the range of $\pm 200$ rad/s while keeping $\nu_x(0)=\nu_y(0)=\dot{\theta}(0)=0$. As $\theta=0.1$ is small, $\cos\theta=0.995\approx 1$, $\sin\theta=0.01\approx0$, so $\dot{\varphi}(0)\in [-200,200]$ affects the value of $\lambda$ as little as $1\%$, so that $\lambda\approx 2\lambda_{\text{thres}}$.

An increase of $\dot{\varphi}(0)$ from 0 to 83.1 shortens the length of the initial synchronisation phase so that $t_{\theta\approx 0}$ goes to zero when $\dot{\varphi}(0)\approx 83.1$. With further increase of $\dot{\varphi}(0)$ the synchronisation phase disappears and the climbing time is shortened.

The initiation time $t_{\theta\approx0}$ also decreases as $\dot{\varphi}(0)$ becomes negative. Remarkably the IC with $\dot{\varphi}(0)\approx 0$ have the longest synchronatisation phase for the whole range $\dot{\varphi}(0)\in[-200,200]$. It may be related to the fact that the total energy is close to its minimal value when $\lambda=2\lambda_{\text{thres}}$ is kept fixed.

\bigskip

\noindent Another distinguished range of angular velocities is $\dot{\varphi}(0)\in[195,200]$ (with $\theta(0)=0.1$) that gives rise to solutions with very small amplitude oscillations for all dynamical variables (see Fig. 8). The increase of $\theta(t)$ becomes monotonous starting from $\theta=0.8$. Solutions with $\dot{\varphi}(0)\in[195,200]$ become even more smooth when $\theta(0)=0.01$, and the increase of $\theta(t)$ is monotononous starting from time $0.3$s and angle $\theta(t)=0.02$. Their initiation time is zero and their ending time is about $t_{\theta\approx\pi}=7.77$.

For testing stability of qualitative behaviour of the reference solution and solutions with $\dot{\varphi}(0)\in[-200,200]$ we have varied the additional initial conditions for $\nu_x(0)$, $\nu_y(0)$ and $\dot{\theta}(0)$ within the maximal range of variability of $\nu_x(t)\in[-0.3,0.3]$, $\nu_y(t)\in[-0.8,0.8]$ and $\dot{\theta}(t)\in[-15,15]$ displayed by the reference solution \eqref{IC_ref}.

The general pattern is that with adding large nonzero $\nu_x(0)$, $\nu_y(0)$ or $\dot{\theta}(0)$ the synchronisation phase may reappear, the climbing time remains below 7-8 s, and the main qualitative features of inverting remain intact.
The price for taking large values of $\nu_x(0)$, $\nu_y(0)$ or $\dot{\theta}(0)$ is that the dynamical variables $\nu_x(t)$, $\nu_y(t)$ and $\dot{\theta}(t)$ oscillate strongly with amplitudes staying within the same maximal range variability $\nu_x(t)\in[-0.3,0.3]$, $\nu_y(t)\in[-0.8,0.8]$ and $\dot{\theta}(t)\in[-15,15]$.
Small perturbations not exceeding $10\%$ of the range of variability of $\nu_x(t)$, $\nu_y(t)$ and $\dot{\theta}(t)$ do not change much the parameters and behaviour of inverting solutions.

Changing of the initial inclination angle $\theta(0)$ affects the amplitude of oscillation for the dynamical variables $\nu_x(t)$, $\nu_y(t)$ and $\dot{\theta}(t)$ and the climbing time. A general rule is that decreasing $\theta(0)$ below $0.1$ reduces oscillations, smooths out solution curves and increases the climbing time. Increase of $\theta(0)$ makes amplitudes of oscillations larger and decreases the climbing time. Beyond the angle $\theta(0)=0.3$ the term $\alpha-\cos\theta$ starts to play a role and $\omega_3(0)$ has to be increased to keep $\lambda\approx 2\lambda_{\text{thres}}$. 

\subsection{Transfer of energy between modes and illustration of the Main equation for the Tippe Top approach}

\begin{figure}[ht]
\begin{center}
\includegraphics[width=\textwidth,height=0.4\textheight]{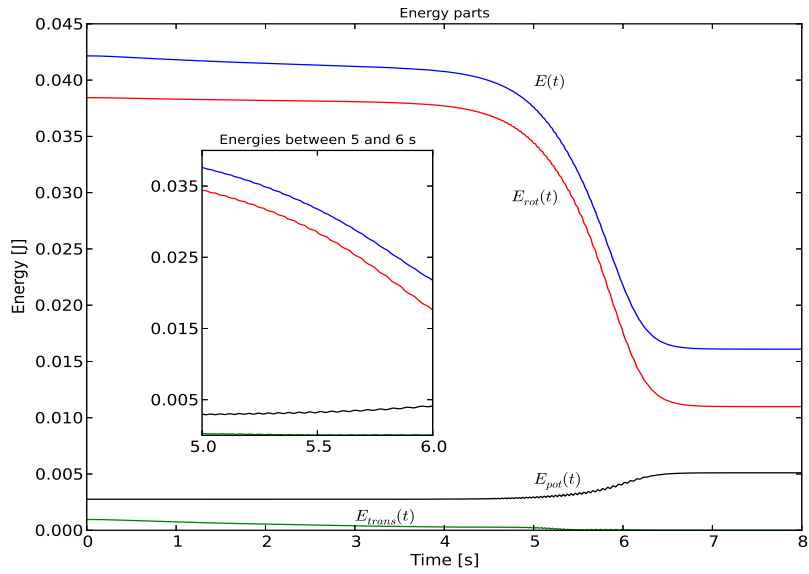}
\caption{Evolution of the energy $E(t)$ and its components $E_{\text{rot}}$, $E_{\text{pot}}$ and $E_{\text{trans}}$. The magnification shows oscillatory behaviour of the energy parts.\label{energy_plot}}
\end{center}
\end{figure}

The total energy consists of three components:
\begin{align}
\label{energy_split}E&=E_{\text{trans}}+E_{\text{rot}}+E_{\text{pot}}=\frac{1}{2}m\mathbf{\dot{s}}^2+\frac{1}{2}\boldsymbol{\omega}\cdot\mathbf{L}+mg\mathbf{s}\cdot\hat{z}\nonumber\\
&=\frac{1}{2}m\left[(\nu_x\cos\theta-R\dot{\theta}(\alpha-\cos\theta))^2+(\nu_y-R\sin\theta(\omega_3+\dot{\varphi}(\alpha-\cos\theta)))^2+(\nu_x\sin\theta+R\dot{\theta}\sin\theta)^2\right]\nonumber\\
&+\frac{1}{2}\left[I_1(\dot{\theta}^2+\dot{\varphi}^2\sin^2\theta)+I_3\omega_{3}^2\right]+mgR(1-\alpha\cos\theta),
\end{align}
since $\mathbf{\dot{s}}=\mathbf{v}_A-\boldsymbol{\omega}\times\mathbf{a}$. Their time-dependence is shown in Fig. \ref{energy_plot} where we see that $E_{\text{trans}}$, $E_{\text{rot}}$ and $E_{\text{pot}}$ add up to the total energy $E(t)$ that is monotonically decreasing.

The components $E_{\text{trans}}$, $E_{\text{rot}}$ and $E_{\text{pot}}$ behave in an oscillatory way since their time derivatives (see below) change sign many times. These oscillations are not directly seen in the main graph because the relative variation of the energy components are small. They become visible in the magnification of the curve (see insert in Fig. \ref{energy_plot}).

In Fig. \ref{energy_plot} the potential energy $E_{\text{pot}}$ (light blue line) increases and the numerical ratio of energies seen in the figure $E_{\text{pot}}(t=8)/E_{\text{pot}}(t=0)=\frac{0.005106}{0.002755}=1.8534$ differs from the theoretical value $\frac{mgR(1+\alpha)}{mgR(1-\alpha)}=\frac{1.3}{0.7}=1.8571$ only by $0.2\%$.
The translational energy $E_{\text{trans}}$ for the chosen IC \eqref{IC_ref} is initially small, about $0.34E_{\text{pot}}(t=0)$ and it goes to zero as TT inverts and approaches asymptotically the inverted spinning state having fixed center of mass. During the inversion the greatest part of the energy is contained in the rotational energy mode $E_{\text{rot}}(t)$.   

The transfer of energy between energy modes becomes visible when we calculate the derivative
\begin{align}
\dot{E}&=\frac{d}{dt}\left(\frac{1}{2}m\mathbf{\dot{s}}^2+\frac{1}{2}\boldsymbol{\omega}\cdot\mathbf{L}+mg\mathbf{s}\cdot\hat{z}\right)\nonumber\\
&=\mathbf{\dot{s}}\cdot(\mathbf{F}_{R}+\mathbf{F}_{f}-mg\hat{z})+\boldsymbol{\omega}\cdot\left[\mathbf{a}\times(\mathbf{F}_{R}+\mathbf{F}_{f})\right]+mg\mathbf{\dot{s}}\cdot\hat{z}\nonumber\\
&=(\mathbf{v}_{A}-\boldsymbol{\omega}\times\mathbf{a})\cdot(\mathbf{F}_{R}+\mathbf{F}_{f})-mg\mathbf{\dot{s}}\cdot\hat{z}+\boldsymbol{\omega}\cdot\left[\mathbf{a}\times(\mathbf{F}_{R}+\mathbf{F}_{f})\right]+mg\mathbf{\dot{s}}\cdot\hat{z}\nonumber\\
&=\mathbf{v_A}\cdot\mathbf{F}_f-(\boldsymbol{\omega}\times\mathbf{a})\cdot(\mathbf{F}_{R}+\mathbf{F}_{f})+\boldsymbol{\omega}\cdot\left[\mathbf{a}\times(\mathbf{F}_{R}+\mathbf{F}_{f})\right]=\mathbf{v_A}\cdot\mathbf{F}_f=-\mu g_n\mathbf{v}_{A}^2,
\end{align} 
by using the dynamical Eqs. \eqref{TTequ}, the equality $\dot{\boldsymbol{\omega}}\cdot\mathbf{L}=\boldsymbol{\omega}\cdot\mathbf{\dot{L}}$ valid due to axial symmetry of TT, and $\mathbf{v}_{A}\cdot\hat{z}=0$. In this calculation $\mathbf{v_A}\cdot\mathbf{F}_f=-\mu g_n\mathbf{v}_{A}^2$ is the rate of frictional loss of energy, $mg\mathbf{\dot{s}}\cdot\hat{z}$ is the rate of energy transfer from the translational part into potential energy $E_{\text{pot}}$. The term $\boldsymbol{\omega}\cdot\left[\mathbf{a}\times(\mathbf{F}_{R}+\mathbf{F}_{f})\right]$  is the work performed in unit time by the torque $\mathbf{a}\times(\mathbf{F}_{R}+\mathbf{F}_{f})$ and it transfers energy between the translational and rotational components.

\begin{figure}[ht]
\begin{center}
\includegraphics[width=\textwidth,height=0.4\textheight]{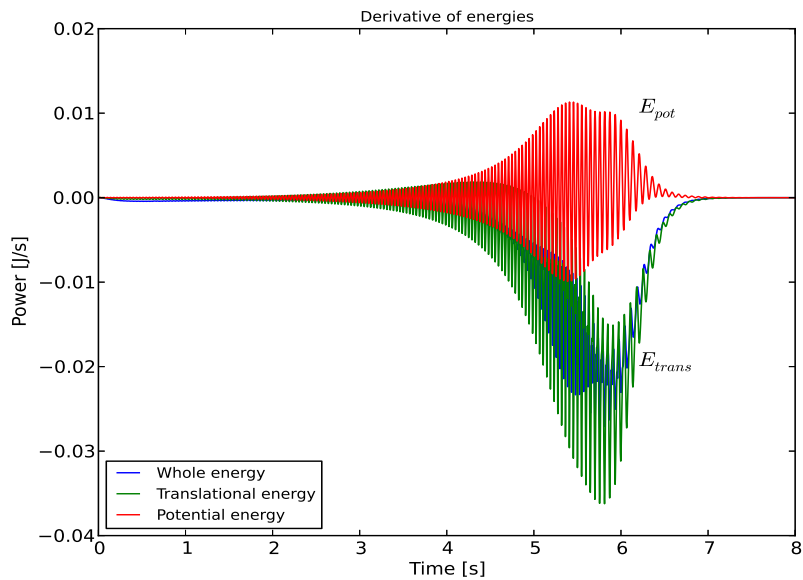}
\caption{Evolution of the derivative of the energy $\dot{E}(t)$ (blue) and the derivative of the translational and potential components $\dot{E}_{\text{trans}}(t)$ (green) and $\dot{E}_{\text{pot}}(t)$ (red). \label{energy_diff_plot}}
\end{center}
\end{figure}

The graphs for $(\dot{E},\dot{E}_{\text{pot}},\dot{E}_{\text{trans}})$ (Fig. \ref{energy_diff_plot}), which describe velocity of energy transfer between modes, are oscillatory. The derivative for the total energy $\dot{E}(t)=-\mu g_n\mathbf{v}_{A}^2\leq 0$ is non-positive but it becomes close to zero at some instants of time when both $\nu_x(t)$ and $\nu_y(t)$ are close to zero in a neighborhood of the initiation time $t_{\theta\approx0}=3.2$ s. The graph of $\dot{E}_{\text{trans}}$ is predominantly negative and has general shape similar to $\dot{E}(t)$. This reflects the fact that most of the energy lost to friction comes from the rotational and the translational components.

\bigskip

\begin{figure}[ht]
\begin{center}
\includegraphics[width=\textwidth,height=0.4\textheight]{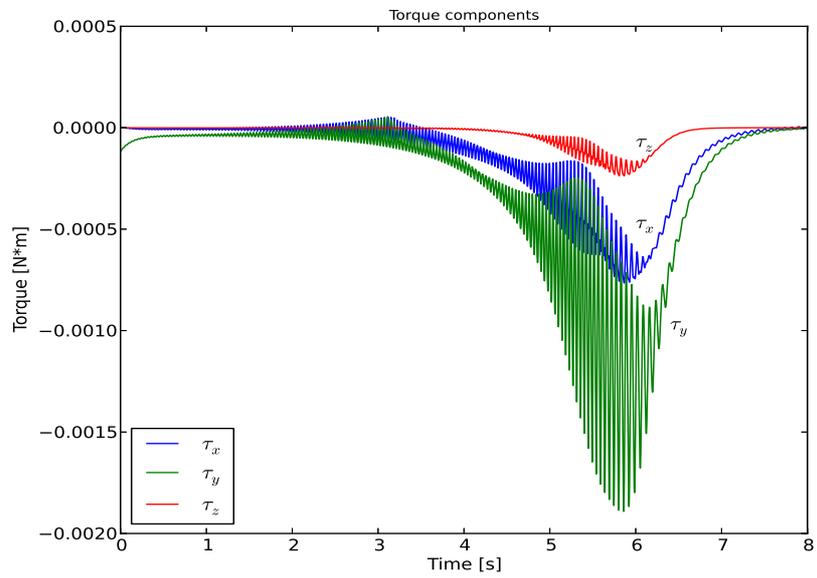}
\caption{The components of the torque $\boldsymbol{\tau}=\mathbf{a}\times(\mathbf{F}_R+\mathbf{F}_f)$. All components are negative almost all the time.\label{torque_plot}}
\end{center}
\end{figure}

\noindent Figure \ref{torque_plot} shows time-dependence of components of the torque vector calculated w.r.t. the center of mass
\begin{equation}
\boldsymbol{\tau}=\mathbf{a}\times(\mathbf{F}_R+\mathbf{F}_f)=\mathbf{a}\times(g_n\hat{z}-\mu g_n\mathbf{v}_A)
\end{equation}
These components are all predominantly negative w.r.t. the chosen moving reference frame $(\hat{x},\hat{y},\hat{z})$. The components $\boldsymbol{\tau}_x=-R(1-\alpha\cos\theta)\mu g_n\nu_y$ and $\boldsymbol{\tau}_z=-R\alpha\mu g_n\nu_y\sin\theta$ are negative whenever $\nu_y>0$, and $\boldsymbol{\tau}_y=-R\alpha g_n\sin\theta+R\mu g_n\nu_x(1-\alpha\cos\theta)$ is negative whenever $\nu_x<0$ or $\nu_x$ is sufficiently small.

\begin{figure}[hbt]
\begin{center}
\includegraphics[width=\textwidth,height=0.6\textheight]{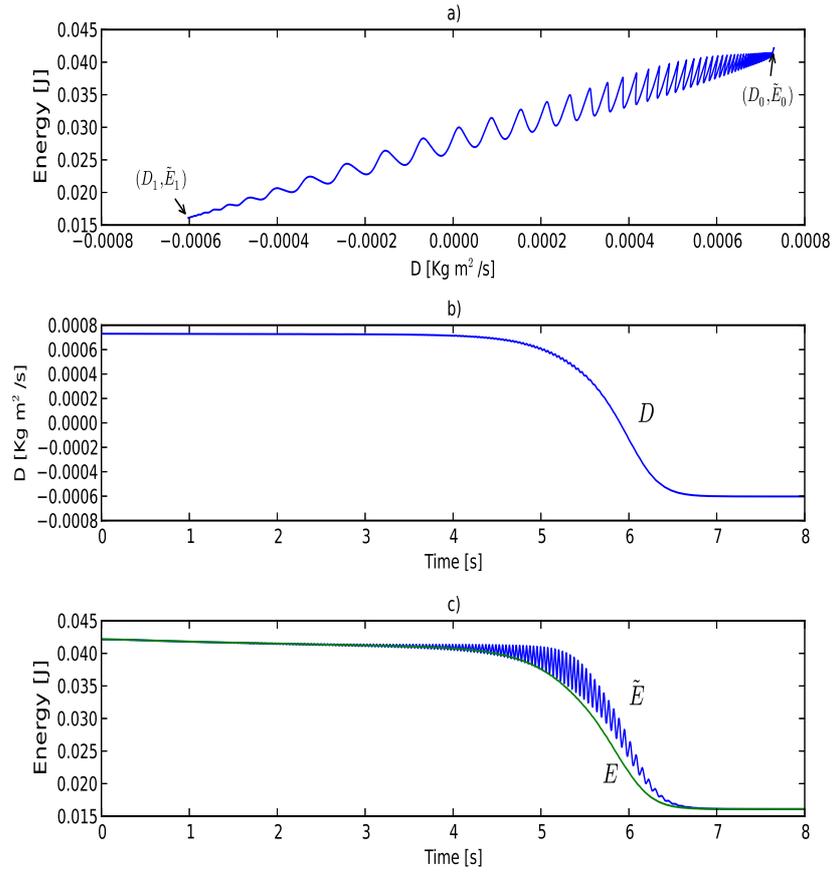}
\caption{Plot $a$ shows the curve $(D(t),\tilde{E}(t))$, calculated by integrating Eqs. \eqref{ddth}-\eqref{nu_y} and Eq. \eqref{equ_for_g_n} with the reference IC used in figure \ref{total_plots1}. Plot $b$ shows the evolution of the function $D(t)$ and plot $c$ shows the energies $E(t)$ (green) and $\tilde{E}(t)$. \label{total_plots3}}
\end{center}
\end{figure}

As explained in section 1 each inverting solution of TT equations satisfies its own Main Equation for the Tippe Top $\tilde{E}(t)=\frac{1}{2}(2g(\cos\theta)\dot{\theta}^2)+V(\cos\theta,D(t),\lambda)$ with $\tilde{E}(t)$, $D(t)$ calculated from \eqref{Etilde_def} and \eqref{Routh_def}. Remember that $\tilde{E}(t)$ is the part of the whole energy \eqref{energy_split} not depending on the gliding velocity $\mathbf{v}_{A}$. The function $D(t)$ decreases from the value $D_0\approx\frac{\lambda\sqrt{d(1)}}{R(1-\alpha)}=9.381\cdot10^{-4}$ to $D_1\approx -\frac{\lambda\sqrt{d(-1)}}{R(1+\alpha)}=-6.667\cdot 10^{-4}$ and in diagram \ref{total_plots3}a the curve $(D(t),\tilde{E}(t))$ goes from the initial boundary value ($D_0\approx 9.381\cdot10^{-4},\tilde{E}_0\approx 4.912\cdot 10^{-2}$) to the final boundary value $(D_1\approx-6.667\cdot 10^{-4},\tilde{E}_1\approx 1.734\cdot10^{-2})$.  The shape of the curve reflects the oscillatory behaviour of $\tilde{E}(t)$ while $D(t)$ (in Fig. \ref{total_plots3}b) is almost monotonously decreasing since $D(t)=I_3\omega_3(t)\sqrt{d(\cos\theta(t))}$.

The blue curve in Fig. \ref{total_plots3}c shows oscillatory behaviour of the modified rotational energy $\tilde{E}(t)$ in relation to the green curve $E(t)$ representing total energy. The modified rotational energy $\tilde{E}(t)$ is consistently larger than $E(t)$ during inversion. This does not contradict the conservation of energy when we look closer at $\tilde{E}=E-\frac{1}{2}m\mathbf{v}_{A}^2+m\mathbf{v}_{A}\cdot(\boldsymbol{\omega}\times\mathbf{a})$. For an inverting TT the angular velocity $\boldsymbol{\omega}$ remains close to be parallel with the $\hat{z}$ axis and the product $\boldsymbol{\omega}\times\mathbf{a}$ points behind the plane of the picture in Fig. \ref{TT_diagram}. The direction of rotation of TT (with $\boldsymbol{\omega}$ almost parallel to $\hat{z}$) causes $\nu_y=\mathbf{v}_{A}\cdot\hat{y}$ to be positive and to point also behind the plane of Fig. \ref{TT_diagram}. 
 
Thus the term
\begin{align}
m\mathbf{v}_{A}\cdot(\boldsymbol{\omega}\times\mathbf{a})-\frac{1}{2}m\mathbf{v}_{A}^2&\approx\frac{1}{2}m(\nu_x\hat{x}+\nu_y\hat{y})\cdot\left(2(\boldsymbol{\omega}\times\mathbf{a})\hat{y}-(\nu_x\hat{x}+\nu_y\hat{y})\right)\nonumber\\
\label{est_v_A}&\approx-\frac{1}{2}m\nu_{x}^2+\frac{1}{2}m\nu_y(2|\boldsymbol{\omega}\times\mathbf{a}|-\nu_y)
\end{align}
may be positive if $\nu_{y}>\nu_x$ and $2|\boldsymbol{\omega}\times\mathbf{a}|-\nu_y=2|\boldsymbol{\omega}||\mathbf{a}|\sin\theta-\nu_y>0$ is sufficiently large.
In Fig. \ref{total_plots2}b we cannot estimate the size of each term, but as we see in Fig.~\ref{TT_diagram}, $|\mathbf{a}|$ is growing during inversion, $|\boldsymbol{\omega}|$ is decreasing and $\sin\theta$ is growing until the axis $\mathbf{\hat{3}}$ passes $\theta=\frac{\pi}{2}$. This may altogether keep the term $(2|\boldsymbol{\omega}\times\mathbf{a}|-\nu_y)$ positive and sufficiently large to make \eqref{est_v_A} positive. This is consistent with the graph \ref{total_plots2}c where the difference $\tilde{E}-E$ becomes largest in the middle of inversion when $\theta\approx\frac{\pi}{2}$.

\section{Summary and conclusions}

Application of the LaSalle principle and stability analysis of asymptotic solutions provides necessary conditions for physical parameters and for initial conditions (IC) so that the Tippe Top inverts. Analysis of dynamical behaviour of inverting solutions is indeed difficult since one deals with a nonlinear, nonintegrable dynamical system of 6 degrees of freedom.

In this paper we have numerically studied properties of solutions starting at small initial inclination angle $\theta(0)=0.1$ rad and with $\lambda\approx 2\lambda_{\text{thres}}$ to learn how such solutions depend on the choice of the remaining IC. Numerical simulations confirm that all such sample solutions invert. Closer analysis of numerical solutions have also shown new interesting features of dynamical behaviour and how they depend on changes of the remaining IC while keeping $\theta(0)=0.1$ and $\lambda\approx 2\lambda_{\text{thres}}$.

\bigskip 

\noindent For solutions with $\nu_x(0)=\nu_y(0)=\dot{\theta}(0)=0$ and $\omega_3(0)=155$ rad/s, a study of dependence on $\dot{\varphi}(0)\in[-200,200]$ shows that solutions with $\dot{\varphi}(0)>83$ rad/s start to invert directly and solutions with $\dot{\varphi}(0)\leq 83$ require a synchronisation time-interval for $\dot{\varphi}(t)$, $\nu_x(t)$ and $\nu_y(t)$ before $\theta(t)$ starts to climb. The climbing {\it initiation time} $t_{\theta\approx0}$ is distinguished:
\begin{itemize}
\item[a)] by the inclination angle being close to zero $\theta(t_{\theta\approx0})\approx0$,
\item[b)] by high amplitude oscillations of $\dot{\varphi}(t)$ in the vicinity of $t_{\theta\approx 0}$,
\item[c)] by reversal of amplitude of oscillations of $\dot{\varphi}(t)$ from large negative values to large positive values of $\dot{\varphi}(t)$,
\item[d)] by the fact that at time $t_{\theta\approx0}$ the graph of $\nu_{y}(t)$ crosses the graph of $\nu_x(t)$ from below and increases by 1-2 orders of magnitude.
\end{itemize}

For all tested IC the climbing time-interval ends at the {\it ending time} $t_{\theta\approx\pi}$ distinguished by another reversal of high amplitude oscillations of $\dot{\varphi}(t)$. These oscillations are well visible in all graphs of a fully inverting Tippe Top. At $t_{\theta\approx\pi}$ the graph of $\nu_y(t)$ crosses the graph of $\nu_x(t)$ from above. So for the climbing solutions there appears to exist a naturally defined {\it climbing time} interval $[0,t_{\theta\approx\pi}-t_{\theta\approx 0}]$ or $[0,t_{\theta\approx\pi}]$ (if $\dot{\varphi}(0)>83$ and the climbing starts immediately) during which the angle $\theta$ increases to $\pi$. Perturbations of the reference IC \eqref{IC_ref} with nonzero values of $\nu_x(0)$, $\nu_y(0)$ and $\dot{\theta}(0)$, belonging to the range of variability of $\nu_x(t)$, $\nu_y(t)$ and $\dot{\theta}(t)$ for the reference solution, preserves the same main features of inverting solutions but amplitudes of oscillations for all dynamical variables usually increases. This may be related to the fact that the total energy for the reference solution is close to the minimal value of total energy among solutions having $\lambda\approx 2\lambda_{\text{thres}}$.

\bigskip

\noindent We have discussed how the observed dynamical behaviour of solutions is related to of the dynamical equations (2)--(6). Such analysis does not provide a proof of the observed dynamical features but it sheds light on the relationships between the main dynamical variables during the inversion. This may serve as a good starting point for further numerical experiments studying dynamics of inverting solutions of a Tippe Top.

We have also illustrated how energy is transfered during the inversion between three energy modes $E_{\text{trans}}(t)$, $E_{\text{rot}}(t)$ and $E_{\text{pot}}(t)$ that together add up to the monotonously decreasing total energy $E(t)$. This transfer is due to the torque that changes the components of angular velocity and performs work needed for energy transfer.

The description of TT inversion through the Main Equation of the Tippe Top has been illustrated by the graph of the modified rotational energy $\tilde{E}(t)$ versus total energy $E(t)$ and by a picture of the curve $(D(t),\tilde{E}(t))$ that controls deformation of the effective potential $V(\cos\theta,D(t),\lambda)$ and motion of $\theta(t)$ inside the potential. 

\bigskip

\noindent The numerical study of inverting solutions confirms the predictions about inversion of TT known from analysis of stability of asymptotic solutions. Additionally it confirms oscillatory behaviour of the angle $\theta(t)$ as predicted by the Main Equation for the Tippe Top approach. It also reveals further intricate features of Tippe Top behaviour that are only partially understood on the basis of dynamical equations. They asks for further research in this direction.   

\bibliography{mybib}

\begin{figure}[ht]
\begin{center}
\includegraphics[width=\textwidth,height=0.6\textheight]{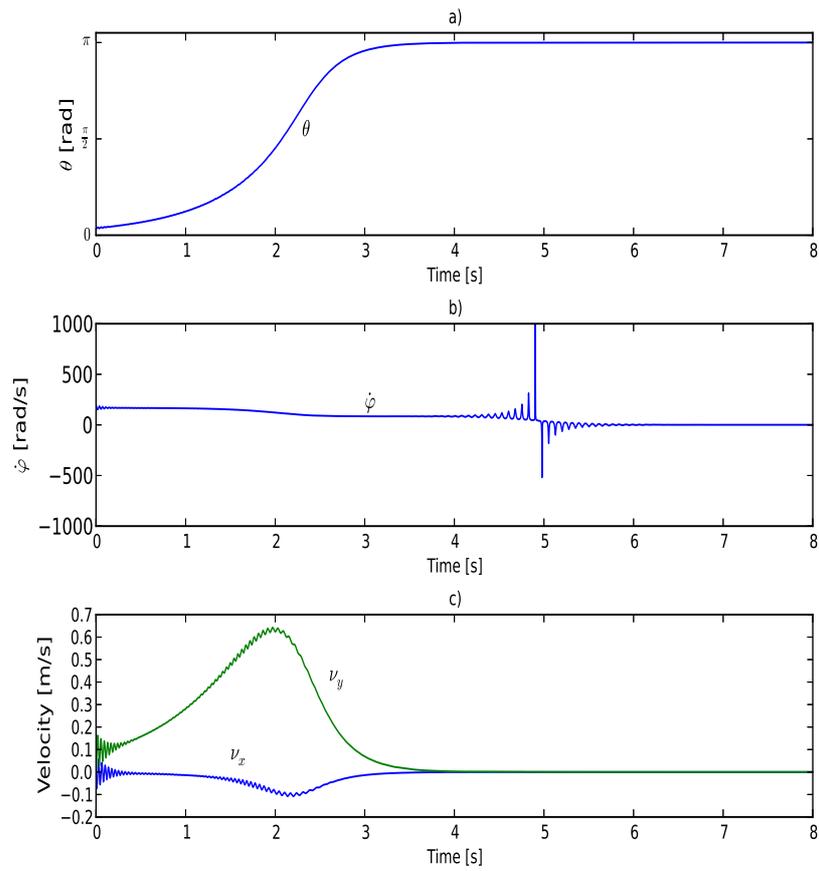}
\caption{Plot for adjusted initial conditions \eqref{IC_ref} with $\dot{\varphi}(0)=195$ rad/s. The behaviour of the dynamical variables is exceptionally smooth. \label{phi_plot}}
\end{center}
\end{figure}

\end{document}